\documentclass[12pt]{article}
\usepackage{amsmath}
\usepackage{amssymb}

\usepackage[T1]{fontenc} 
\usepackage{sets pace}
\usepackage{comment}
\usepackage{colortbl}
\definecolor{text1}{cmyk}{1,.65,0,0} % blue text colour
\definecolor{text2}{rgb}{1,0,0} % red text colour
\definecolor{text3}{cmyk}{0,0,0,1} % black text colour
\definecolor{text4}{cmyk}{0,0,0,0.5} % grey text colour
\definecolor{text5}{cmyk}{1.0,0.0,1.0,0} % green text colour

\usepackage{amsfonts}
\usepackage{amsthm}
\usepackage{subcaption}
\usepackage{graphicx}
\usepackage{caption}
\usepackage{subcaption}
\usepackage{latexsym}
\usepackage[most]{tcolorbox}
\makeatletter
\renewcommand{\@seccntformat}[1]
{\csname the#1\endcsname.\enspace}
\makeatother
\newtheorem{definition}{Definition}
\newtheorem{theorem}{Theorem}
\newtheorem{lemma}{Lemma}
\newtheorem{remark}{Remark}
\newtheorem{corollary}{Corollary}

\begin{document}

\large
\begin{center}
{\bf  On the Bayesian analysis of a non-identifiable binomial model \footnote{ \today. }}   
\end{center}
\normalsize
%begin{tabular}
\begin{center}
{\sc \'Eric Marchand } \\
{\it   Universit\'e de
    Sherbrooke, D\'epartement de math\'ematiques, Sherbrooke Qc,
    CANADA, J1K 2R1 \quad (e-mail: eric.marchand@usherbrooke.ca) } \\
    \vspace*{0.3cm}
    {\sc M.C. Jones } \\
\textit{School of Mathematics and Statistics, The Open 
University, Milton Keynes MK7 6AA, UK} \\
\end{center}

\vspace*{-0.5cm}
\begin{abstract}
\noindent
We provide analysis for the posterior distribution and expectation of $(p_1, p_2)$ where $Y|p_1,p_2 \sim \hbox{Binomial}(n, p_1 p_2)$ and $ (p_1, p_2)$ is uniformly distributed on the unit square $[0,1]^2$.   We exhibit: \textbf{(i)} various mixture representations for the posterior distributions of $p_1|p_2,y$ and $p_2|y$ which involve truncated Beta, truncated binomial and truncated negative binomial distributions, \textbf{(ii)} a simple and attractive ``harmonic number distribution'' for the marginal pmf $p(y)$ of $Y$, \textbf{(iii)} large sample size $n$ approximations for $p(y)$, the posterior expectations $\mathbb{E}(p_i|y)$ and variances $\mathbb{V}(p_i|y)$, as well as for the normalization constant in the joint posterior density, and \textbf{(iv)} exact and large sample approximations for the marginal posterior density $p_2|y$ indicative of the absence of posterior contraction. 
\end{abstract}

\vspace{0.1cm} 


\section{Introduction}

We consider the Bayesian model:
\begin{equation}
\label{model}
Y|p_1,p_2 \sim \hbox{Binomial}(n, p_1 p_2) \; \hbox{with } \;  p_1, p_2 \sim^{\!\!\!\!\!\!\!\! i.i.d.} \mathcal{U}(0,1)\,, 
\end{equation}
and extract various interesting representations relative to the Bayesian posterior distribution and expectation.  The above binomial problem arises naturally in a basic setting with independent, but otherwise arbitrary, events $A$ and $B$ such that $\mathbb{P}(A)\,=\, p_1$ and $\mathbb{P}(B)\,=\, p_2$, and where $n$ independent replicates of the indicator function $\mathbb{I}_{A \cap B}$ are observed.   In reliability, the outcome $Y$ represents the success of a two-component series system, and prior/posterior information for $(p_1,p_2)$ relates to the prior/posterior reliability of subsystems.   For Bayesian analysis of more complex series/parallel systems and their components, we refer to Martz et al. (1988), as well as to Martz \& Walker (1990), as well as references therein. 

The inclusion of the uniform prior density in  (\ref{model}) appears in a recent paper (Surjanovic et al., 2025), where the authors introduce and illustrate the use of the software package Pigeons.jl,  which is designed to sample from difficult or intractable distributions including Bayesian posterior distributions. The authors refer to the above as a non-identifiable model for a  ``coinflip'' data set since the effects of the parameters $p_1$ and $p_2$ cannot be distinguished.  From an inferential point of view, the non-identifiability of the model leads to an important challenges, in particular with regards to implications for Bayesian analysis.  An example of a discussion of some of the issues that arise (e.g., choice of prior and its sensitivity), and moreover illustrated for a similar model $Y|p_1,p_2 \sim   
\hbox{Binomial}(n, p_1 + p_2)$ with $p_1 \geq 0\,, p_2 \geq 0\,, \,p_1+p_2 \leq 1$, is given by  Neath \& Samaniego (1997).
  
We report here on interesting features which include: 
\begin{enumerate}
\item[ (i)]   an exact and simple form for the marginal probability and cumulative distributions (cdf) functions of $Y$ and related normalization constant, as well as large sample approximations;
\item[ (ii)]     Beta mixture representations for the marginal posterior density $\pi(p_2|y)$ as well as for the conditional posterior distribution $\pi(p_2|p_1,y)$, the latter also relating to a truncated Beta distribution which we introduce;
\item[ (iii)]   an exact and simple expression for the posterior expectation $\mathbb{E}(p_2|y)$, as well as large sample approximations;
\item[ (iv)]    an asymptotic expression for the posterior variance $\mathbb{V}(p_2|y)$ for 
 $y=\lfloor n \alpha \rfloor$, $\alpha \in (0,1)$ and $n \to \infty$, shown to be strictly positive and indicative of the absence of posterior contraction in such a setting;
 \item[ (v)]  an alternate expression for the posterior density $\pi(p_2|y)$ which involves the cdf of a Beta distribution, as well as an approximate large sample version which involves the cdf of a normal distribution illustrating again the absence of contraction;
\item[ (vi)]  an intimate connection with harmonic numbers in representations (i), (iii) and (v) above.
\end{enumerate}    

\noindent  In relationship to the points above, the manuscript is organized as follows.  Sections 2 and 3 present findings for (i) and (ii), respectively, while Section 4 relates to (iii).  Section 5 has findings associated with (i), (iv) and (v).
 
Although model (1.1) is a simple one and well-formulated for computational evaluations, we believe that the precise posterior and marginal distribution findings presented in this work are particularly welcomed, quite revealing of the underlying finite sample and asymptotic structures, and potentially useful for the validation and implementation of numerical algorithms.

\section{Marginal distribution of $Y$ and normalization constant}

Clearly, the variables $p_1$ and $p_2$ are interchangeable both for the binomial model (\ref{model}) and the prior density, so that for a given observation $y$ the joint posterior density will be exchangeable and the resulting marginal posterior distributions will necessarily match. 
  A collapsed version of the model is $Y|w \sim \hbox{Binomial}(n, w)$ with $W$ distributed as the product $p_1 p_2$.   This has density $f(w) \, = \, - \ln (w) \, \mathbb{I}_{(0,1)}(w)$ and moments $\mathbb{E}(W^k) \, = \, \frac{1}{(k+1)^2}, k \geq 0$.  
With $\mathbb{E}(Y|w)\,=\, nw$ and $\mathbb{V}(Y|w)\,=\, nw(1-w)$, one obtains the marginal expectation and variance $\mathbb{E}(Y)\,=\, \frac{n}{4}$ and $\mathbb{V}(Y) \, = \, \frac{7 n^2}{144} + \frac{5n}{36}$.  

\noindent We proceed with representations for the marginal probability mass function (pmf) and cdf of $Y$.  As shown below, these are intimately related to harmonic numbers and summation.   We denote $H_m$ as the $m^{th}$ harmonic number (with the convention that $H_0=0$) given by $H_j\,=\, \sum_{k=1}^{m} \frac{1}{k}$ for $m \in \mathbb{N}_+$.  We will make use of the following property:  
\begin{equation}
\label{harmonicsummation}
  \sum_{k=0}^y H_k\,=\, (y+1) \, (H_{y+1}-1) \,, \hbox{ for } y \in \mathbb{N}\,,
\end{equation}
as well as the bounds and asymptotic evaluations
\begin{align}
\label{harmonicbounds}
& \ln(m) \, + \, \gamma \, + \, \frac{1}{2(m+1)} \leq H_m \leq \ln(m) \, + \, \gamma \, + \, \frac{1}{2 m} \,, \\ 
\label{harmoniclimits}& \lim_{m \to \infty} \big\{H_m - \ln (m)\big\} = \gamma = 0.577216.. \hbox{ being Euler's constant}, \, \\
\label{harmonicapprox}
& \hbox{and } \, H_m \simeq \ln m\,+\,\gamma+\frac1{2m}\,, \hbox{ or } H_m \simeq \ln m\,+ \,\gamma+\frac1{2m}-\frac1{12m^2}\,.
\end{align}
%\begin{equation} \label{harmonicapprox}
%H_n \simeq \ln n\,+\,\gamma+\frac1{2n}\,, \hbox{ or } H_n \simeq \ln n\,+ \,\gamma+%\frac1{2n}-\frac1{12n^2}\,.
%\end{equation}
These two approximations may be obtained from Euler's alternating series expansion and are very accurate, with errors bounded in absolute value by $1/12m^2$ and $1/120m^4$, respectively.

\begin{lemma}
\label{lemmamarginalYchris}
The marginal pmf and cdf of $Y$ are given by
\begin{align}
\label{marginalpmf}
p(y) \, & = \, \frac{H_{n+1}-H_y}{n+1}\,\,{\mathbb{I}}_{\{0,1,...,n\}}(y)\,, \\
\label{marginalcdf} \mathbb{P}(Y \leq y) \, & = \,  \frac{y+1}{n+1} \, (H_{n+1}-H_{y+1}+1) \,, \hbox{ for } y \in \{0,1,\dots, n\}\,.  
\end{align}  
%where $H_j$ is the $jth$ harmonic number (with the convention that $H_0=0$) given by 
%$H_j\,=\, \sum_{k=1}^j  \frac{1}{k}$ for $k \in \mathbb{N}_+$.
\end{lemma}
{\bf Proof.}  The second of these expressions follows from (\ref{marginalpmf}) and (\ref{harmonicsummation}) as
\begin{equation}
\nonumber  \sum_{k=0}^y p(k) \,=\, \sum_{k=0}^y  \frac{H_{n+1}-H_y}{n+1}\, = \, \frac{y+1}{n+1} \, H_{n+1} \, - \, \frac{1}{n+1} \, (y+1) \, (H_{y+1}-1)\,,
\end{equation}
which is (\ref{marginalcdf}).  For the pmf, we have
\begin{eqnarray*}
p(y) \, & = & \, - \, \int_0^1  \binom{n}{y} \, w^y (1-w)^{n-y} \, \ln(w) \, dw \\
 &=& -{n \choose y} \frac{d}{da} \left. \int_0^1 w^{y+a}(1-w)^{n-y}\,dw \,\right|_{a=0}\\
 &=& -{n \choose y} \frac{d}{da} \left. \hbox{B}(y+1+a,n-y+1)\, \right|_{a=0}\\
 &=& -\frac{n!}{y!}  \frac{d}{da} \left. \frac{\Gamma(y+1+a)}{\Gamma(n+2+a)}\, \right|_{a=0}\\
 &=& -\frac{n!}{y!} \left\{ \frac{\psi(y+1)y!}{(n+1)!}-\frac{\psi(n+2)y!}{(n+1)!}\right\}\\
 &=& \frac{\psi(n+2)-\psi(y+1)}{n+1}.
\end{eqnarray*}
Hereafter, $\hbox{B}(\cdot, \cdot)$ is the usual Beta function and $\psi(\cdot) = \Gamma'(\cdot)/\Gamma(\cdot)$ is the digamma function. The result follow immediately since, at integer values $j$, $\psi(j+1) = H_{j}-\gamma$.  \qed 

These are simple and attractive forms for the pmf and cdf, otherwise not widely known we believe.  Its pmf is immediately seen to be decreasing, from $p(0) = H_{n+1}/(n+1)$ to $p(n) = 1/(n+1)^2$. A corresponding expression for a series system with $l \geq 2$ independent components is provided in an Appendix.
 
\section{Posterior distribution of $(p_1, p_2)$}

We pursue various properties and decompositions of the joint posterior density of $(p_1, p_2)$ given by
\begin{equation}
\label{jointposteriordensity}
\pi(p_1,p_2|y) \, = \, \frac{1}{p(y)} \binom{n}{y} \, (p_1 p_2)^y \, (1- p_1p_2)^{n-y} \, \, \mathbb{I}_{(0,1)^2}(p_1,p_2)\,. \\
\end{equation}
We first define a `truncated' Beta distribution and present some related mixture distribution representations.
\begin{definition}
\label{definition}
A `truncated'  Beta distribution with parameters $a>0, b>0, \lambda \in [0,1]$, denoted 
$Z \sim \mathbb{B}(a,b,\lambda)$, has density
\begin{equation}
\label{densityweightedbeta}
f(z) \, = \,   \frac{1}{C(a,b,\lambda)} \, z^{a-1} (1-\lambda z)^{b-1} \, \mathbb{I}_{(0,1)}(z)\,.
\end{equation}
\end{definition}
The above family of distributions include standard Beta distributions $\mathbb{B}(a,b,1)$ and  
are related through truncation and multiplication to standard Beta distributions as 
\begin{equation}
\label{representationconditional+truncation}
\frac{Z}{\lambda}| Z \leq \lambda \, =^d \, \frac{Z}{\lambda} \, \mathbb{I}_{(0,\lambda)}(Z)\,| Z \leq \lambda \, =^d \mathbb{B}(a,b,\lambda) \,,
\end{equation}
for $Z \sim \mathbb{B}(a,b,1)$ and $\lambda>0$.

\begin{remark}
\label{GaussHypergeometricBeta}
Densities (\ref{densityweightedbeta}) are contained within the larger family of Gauss hypergeometric densities introduced by Armero \& Bayarri (1994 a,b) in the context of analyzing a queueing system, and which are of the form
\begin{equation}
\nonumber
\label{densityGaussHypergeometricBeta}
f(z) \, = \,   \frac{1}{K} \, \frac{z^{a-1} (1-z)^{c-1}}{(1-\lambda z)^{1-b}} \, \mathbb{I}_{(0,1)}(z)\,,
\end{equation}
for $a>0, b \in \mathbb{R}, c>0, \lambda < 1$, and  with $K= \hbox{B}(a,c) \, _2F_1(1-b, a, a+c, \lambda)$.    Cases with $c=1$ reduce to the densities in (\ref{densityweightedbeta}) with $K=C(a,b,\lambda)\,=\,   
\frac{1}{a} \, _2F_1(1-b, a, a+1, \lambda)$.   For further analysis of the Gauss hypergeometric distribution, we refer to the recent work of Nadarajah \& Kebe (2026), and the references therein. 
\end{remark}
When $b \in \mathbb{N}_+$, so that $1-b$ is a non-positive integer, a usual series representation of the hypergeometric function $_2F_1$ leads to the finite sum
$C(a,b,\lambda) \, = \ \sum_{k=0}^{b-1} \binom{b-1}{k}\frac{(-\lambda)^k}{a+k}$.  
   Here are a couple of different and lesser-known expressions for Definition \ref{definition}'s  normalization constant $C(a,b,\lambda)$ which will lead to a 
mixture representation below in Theorem \ref{theorem}.   
 
\begin{lemma}  We have for $a>0, b \in  \mathbb{N}_+$, and $\lambda \in [0,1]$:
\label{lemmaC}
\begin{equation}
\label{C1}
C(a,b,\lambda) \, = \, \Gamma(a) \, \Gamma(b) \, \sum_{k=0}^{b-1}  \, \frac{\lambda^{b-1-k} \, (1-\lambda)^k }{k! \, \Gamma(a+b-k)} \,,
\end{equation}
and
\begin{equation}
\label{C2}
C(a,b,\lambda) \, = \, \Gamma(a) \, \Gamma(b) \, \sum_{k=0}^{b-1}  \, \frac{\lambda^{k} \, (1-\lambda)^{b-1-k} }{\Gamma(b-k) \,\, \Gamma(a+k+1)}\,.
\end{equation}
\end{lemma}
{\bf Proof.}   Equation (\ref{C2}) follows from (\ref{C1}) with the change of variables $k \to b-1-k$ in the sum.   For (\ref{C1}), set $\rho=\frac{1-\lambda}{\lambda}$.   A binomial expansion, an interchanging of integral and sum, and a Beta integral yield:
\begin{eqnarray*}
C(a,b,\lambda) \, & = & \,  \int_0^1  z^{a-1} (1-\lambda z)^{b-1} \, dz \\
 \, & = & \, \int_0^1  z^{a-1} \, \lambda ^{b-1} \, (1-z+\rho)^{b-1} \, dz \\
\, & = & \, \int_0^1  z^{a-1} \, \lambda ^{b-1} \, \sum_{k=0}^{b-1} \binom{b-1}{k}\, \rho^k \, (1-z)^{b-1-k} \,  dz \\
\, & = & \,  \lambda^{b-1} \, \sum_{k=0}^{b-1} \rho^k \, \binom{b-1}{k} \, \frac{\Gamma(a) \, \Gamma(b-k)}{\Gamma(a+b-k)}\,,
\end{eqnarray*}
which is equal to (\ref{C1}).  \qed

We pursue two useful and novel mixture distribution representations for $Z \sim \mathbb{B}(a,b,\lambda)$ distributions with $b \in \mathbb{N}_+$. They involve
the pmf's 
\begin{equation}
\label{qk's}
q_k \propto  \frac{\rho^k}{k! \, \Gamma(a+b-k)} \, \mathbb{I}_{\{0,1, \dots, b-1\}}(k) \,,
\end{equation}
and 
\begin{equation}
\label{rk's}
r_j \propto  \frac{\Gamma(a+j) \, (1-\lambda)^j}{j!} \, \mathbb{I}_{\{0,1, \dots, b-1\}}(j) \,,
\end{equation}
for parameters $\rho>0$ and $\lambda \in (0,1)$.   The latter case corresponds for all $a>0$ to the pmf of a truncated at $b-1$ negative binomial distribution which we denote $\hbox{TNB}(a, \lambda;b-1)$, while the former pmf corresponds for $a \in \mathbb{N}_+$ to that of a truncated at $b-1$ binomial distribution (i.e., $\hbox{Bin}(a+b-1,\frac{\rho}{1+\rho})$)  which we denote as $\hbox{TBin}(a+b-1, \frac{\rho}{1+\rho};b-1) $.

\begin{lemma}
\label{lemmamixtures}
Consider $Z \sim \mathbb{B}(a,b,\lambda)$ with $b \in \mathbb{N}_+$.    Then, the distribution of $Z$ is representable as the mixtures:
\begin{enumerate}
\item[ (a)] $Z|K=k \sim \mathbb{B}(a,b-k,1)$ with $K \sim \hbox{TBin}(a+b-1, 1-\lambda;b-1)$ for  $a \in \mathbb{N}_+$;
\item[ (b)]  $Z|J=j \sim \mathbb{B}(a+j,b-j,1)$ with $J \sim \hbox{TNB}(a, \lambda;b-1)$ for  $a >0$.
\end{enumerate}
\end{lemma}
{\bf Proof.}
(a) With $\rho=\frac{1-\lambda}{\lambda}$, proceeding as in the proof of Lemma \ref{lemmaC}, we have  
\begin{eqnarray*}
f_{a,b,\lambda}(z)  \, & \propto & \,  z^{a-1} (1-\lambda z)^{b-1} \, \mathbb{I}_{(0,1)}(z) \\
%\, & \propto & \,  z^{a-1} \, \lambda ^{b-1} \, (1-z+\rho)^{b-1} \, \mathbb{I}_{(0,1)}(z) \\
\, & \propto & \,  z^{a-1} \, \sum_{k=0}^{b-1} \binom{b-1}{k}\, \rho^k \, (1-z)^{b-1-k} \, \mathbb{I}_{(0,1)}(z)  \\
\, & \propto & \, \sum_{k=0}^{b-1} q_k  \,f_{a,b-k,1}(z)\,,
\end{eqnarray*}
with the $q_k$'s as in (\ref{qk's}) and the densities $f_{a,b-k,1}$ being those of Beta$(a,b-k)$ distributions.

(b)  Alternatively with a different binomial expansion, the truncated beta density (\ref{densityweightedbeta}) can be written as 
\begin{eqnarray*}
f_{a,b,\lambda}(z)  \, & \propto & \,  z^{a-1}\big\{(1-\lambda)z+(1-z)\big\}^{b-1} \, \mathbb{I}_{(0,1)}(z) \\
\, & \propto & \,  z^{a-1} \sum_{j=0}^{b-1} {b-1 \choose j}(1-\lambda)^{j}z^j(1-z)^{b-1-j} \, \mathbb{I}_{(0,1)}(z) \\
\, & \propto & \, \;\;  \sum_{j=0}^{b-1} r_j  \,f_{a+j,b-j,1}(z)\,,
\end{eqnarray*}
with the $r_j$'s as in (\ref{rk's}), and the densities $f_{a+j,b-j,1}$ being those of Beta$(a+j,b+j)$ distributions. \qed
 
\noindent  We now continue with the main findings of this section.
 
\begin{theorem}
\label{theorem}
Consider $Y|p_1,p_2 \sim \hbox{Binomial}(n, p_1 p_2)$ and uniform prior for $(p_1,p_2)$ on $[0,1]^2$.  Then, we have the following.
\begin{enumerate}
\item[ (a)]  The posterior distribution of $(p_1,p_2)$ is exchangeable and admits the decomposition: 
\begin{equation}
\label{equationbetamixture}
(i) \;  p_1|p_2,y \sim \mathbb{B}(y+1, n-y+1, p_2) \, , \,
(ii)    \hbox{ with } \, \pi(p_2|y) \, = \, \sum_{k=0}^{n-y}  w_k \, f_{y+k+1, n-y-k+1,1}(p_2)\,,
\end{equation}
$f_{a,b}$ being a $\mathbb{B}(a,b,1)$ density, and probabilities 
\begin{equation}w_k \, = \frac{1}{H_{n+1}-H_y}\,  (k+y+1)^{-1} \, \mathbb{I}_{\{0, \ldots, n-y\}}(k) \,.
\end{equation}
%\begin{equation}
%\nonumber  w_k \, = \frac{1}{H_{n+1}-H_y}\,  (k+y+1)^{-1} \, \mathbb{I}_{\{0, \ldots, n-y\}}(k) \,.\end{equation}
 
\item[ (b)] The conditional posterior distribution of $ Z \sim \mathbb{B}(y+1,n-y+1,p_2)$ admits the mixture decomposition $Z|K=k \sim {\rm Beta}(y+1,n-y-k+1,1)$ with $K\sim {\rm TBin}(n+1,1-p_2;n-y)$.  Equivalently,
$$\pi(p_1|p_2,y) = \sum_{k=0}^{n-y} q'_k f_{y+1,n-y-k+1,1}(p_1)$$
where 
$$q'_k \propto  \frac{\{(1-p_2)/p_2\}^k}{k!\,\Gamma(n+1-k)} \, \mathbb{I}_{\{0,1, \dots,n-y\}}(k).$$
\item[ (c)] The conditional posterior distribution of $ Z \sim \mathbb{B}(y+1,n-y+1,p_2)$ admits the mixture decomposition $Z|J=j \sim {\rm Beta}(y+1+j,n-y+1-j)$ with $J\sim {\rm TNB}(y+1,p_2;n-y)$. Equivalently,
$$\pi(p_1|p_2,y) = \sum_{j=0}^{n-y} r'_j \, f_{y+1+j,n-y+1-j,1}(p_1)$$
where 
$$r'_j \propto  {y+j \choose j} (1-p_2)^j \, \mathbb{I}_{\{0,\dots,n-y\}}(j).$$

\item[ (d)]  The marginal posterior density in (\ref{equationbetamixture}) is alternatively  expressible as \\  $\pi(p_2|y)\,=\, \sum_{k'=0}^{n-y} \nu_{k'} f_{y+1,n-y-k'+1,1}(p_2)$, with $$\nu_{k'}\,=\,  w_{n-y-k'} =  \, \,\frac{1}{H_{n+1}-H_y} \, (n-k'+1)^{-1}\mathbb{I}_{\{0,...,n-y\}}(k')\,,$$ i.e., the distribution of $p_2|y$ admits the mixture representation
\begin{equation}
\nonumber
p_2|y,K'=k' \sim {\rm Beta}(y+1,n-y-k'+1,1) \hbox{ with } \mathbb{P}(K'=k')\,= \,  \nu_k'\,. 
\end{equation}
\end{enumerate}
\end{theorem}
{\bf Proof.}  Parts (b) and (c) are immediate applications of Lemma \ref{lemmamixtures} based on the truncated Beta representation in part (a) for the density of $p_1|p_2,y$.  For part (a), we obtain from (\ref{jointposteriordensity})
\begin{equation}
\nonumber
%\pi(p_1,p_2|y) \, & \propto & \,  \binom{n}{y} \, (p_1 p_2)^y \, (1- p_1p_2)^{n-y} \, \, %%\mathbb{I}_{(0,1)^2}(p_1,p_2) \\
 \pi(p_1,p_2|y) \, \propto \, \frac{p_1^y \, (1- p_1 p_2)^{n-y}}{C(y+1, n-y+1, p_2)} \,\, p_2^y  \, C(y+1, n-y+1, p_2) \,,
\end{equation}
which implies (i).   Moreover, we have from (\ref{C2})
\begin{eqnarray*}
\pi(p_2 |y)  \, & \propto \, &  p_2^y  \, \, C(y+1, n-y+1, p_2) \\
\, & \propto \, & \sum_{k=0}^{n-y}  \frac{\Gamma(n+2) \, p_2^{y+k} \, (1-p_2)^{n-y-k}}{\Gamma(n-y-k+1) \, \Gamma(y+k+1)} \,  \	\times \frac{1}{y+k+1} \,,
\end{eqnarray*}
which yields part (ii). 

\noindent Finally, for part (d), we have from the proof of Lemma \ref{lemmamixtures} and (\ref{rk's}) that
\begin{equation}
\nonumber  C(a,b,\lambda) \propto \sum_{k=0}^{b-1}r_k \propto \sum_{k=0}^{b-1} \frac{\Gamma(a+k)(1-\lambda)^k}{k!}\,.
\end{equation}
%\begin{eqnarray*}
%C(a,b,\lambda) &=& \int_0^1 z^{a-1}(1-\lambda z)^{b-1}\,dx\nonumber\\
%&=& \int_0^1 z^{a-1}\sum_{k=0}^{b-1} {b-1 \choose k}(1-\lambda)^{k}z^k(1-z)^{b-1-k}\,dz\nonumber\\
%&=& \sum_{k=0}^{b-1} {b-1 \choose k}(1-\lambda)^{k}\, \frac{\Gamma(a+k)\,\Gamma(b-k)}{\Gamma{(a+b)}}.
%\end{eqnarray*}
Thus, it follows that
\begin{eqnarray*}
\pi(p_2|y) &\propto & p_2^y\, C(y+1,n-y+1,p_2)\\
%&\propto & p_2^y\sum_{k=0}^{n-y} {n-y \choose k}(1-p_2)^{k}\, {(y+k)!\,(n-y-k)!}\\
&\propto & \sum_{k=0}^{n-y} p_2^y (1-p_2)^{k}\, \frac{(y+k)!}{k!}\\
&=& \sum_{k=0}^{n-y} p_2^y (1-p_2)^{n-y-k}\, \frac{(n-k)!}{(n-y-k)!}\\
&\propto & \sum_{k=0}^{n-y} \nu_k f_{y+1,n-y-k+1,1}(p_2)\,. \,\;\;\; \qed
\end{eqnarray*}

Parts (a) and (d) of Theorem \ref{theorem} establish two intriguing Beta mixture representations for the marginal posterior distributions of $p_2$ (or $p_1$).  They involve discrete mixing variables taking on $n-y+1$ values with probabilities proportional to consecutive elements of a harmonic sequence.  Parts (b) and (c) of Theorem \ref{theorem} provide two distinct Beta mixture representations for the conditional posterior distribution of $p_1|p_2,y$ (or $p_2|p_1,y$), which bring into play truncated binomial and truncated negative binomial mixing.   The equivalence of these representations is noticeable.

\noindent  We conclude this section with a couple of additional remarks.  

\begin{remark} 
\label{remarksimulation}
Part (a) of Theorem \ref{theorem} suggests various options for simulating from the posterior distribution of $(p_1,p_2)$.  On one hand, the full conditionals are available for the Gibbs sampler, i.e., (i) and its twin version $p_2|p_1,y \sim \mathbb{B}(y+1, n-y+1, p_1) \,$ which requires simulating from a truncated Beta distribution, which in turn can be done through its cdf, via representation (\ref{representationconditional+truncation}) or as well as with mixture representations in parts (b) or (c) of Theorem \ref{theorem}.   Alternatively, parts (i) and (ii) (or (d)) combined generates exactly one observation under the posterior, with generators of a truncated Beta distributed variable, discrete variable $K$ (or $K'$) with corresponding generated $\hbox{Beta}(k+y+1, n-y-k+1)$ (or $\hbox{Beta}(y+1, n-y-k+1)$)  variable.    
\end{remark}

\begin{remark} 
\label{remarkY=n}
We point out that the posterior distribution of $(p_1,p_2)$ for $Y=n$ simplifies in (\ref{equationbetamixture}) to independent $\mathbb{B}(n+1,1)$ distributed components.  This is somewhat reminiscent of Laplace's celebrated problem of succession where $Y|U \sim \hbox{Binomial}(n,U)$ with $U \sim \mathcal{U}(0,1)$ also leads to a $\hbox{Beta}(n+1,1)$ posterior distribution for $U$ when $Y=n$.  It is also interesting to observe that the independently distributed property of $p_1$ and $p_2$ is a manifestation of a more general property, specific to the posterior distribution for $Y=n$, where $Y|p_1, \dots, p_k \sim \hbox{Binomial}\big(n, \prod_{i=1}^k p_i\big)$ and $p_1, \dots, p_k$ has prior density of the form $\prod_{i=1}^k g_i(p_i)$.  Indeed, for such a situation, one obtains directly
the joint posterior density
\begin{equation}
\nonumber
\pi(p_1,\dots, p_n|n) \, \propto   \, (\prod_{i=1}^k p_i)^n \, \prod_{i=1}^k g_i(p_i) \,, 
\end{equation} 
which implies independently distributed components $p_i$ with densities
$\pi(p_i|n) \, \propto   \, p_i^n \, g_i(p_i) \, \hbox{ for } i=1,\dots, k\,.$
\end{remark}

\section{Posterior expectation}

We now turn our attention towards the posterior expectation $\mathbb{E}(p_2|y)$ with  the following representation bringing into play a harmonic mean and leading to a simple large $n$ approximation.
\begin{corollary}
\label{corollaryposteriorexpectation}
Consider model (\ref{model}) and suppose that $Y=y$ has been observed.
\begin{enumerate}
\item[ (a)]   The marginal posterior expectation of $p_2$ (or $p_1$) is given by
\begin{equation}
\label{expectationharmonic1}   \mathbb{E}(p_2|y) \, = \, \frac{(n-y+1)/(n+2)}{H_{n+1} - H_y},
\end{equation}
or, alternatively,
\begin{equation}
\label{expectationharmonic2}
\mathbb{E}(p_2|y) \, = \, \frac{1}{n+2} \; \mu_H\{y+1, y+2, \ldots, n+1  \}\,,
\end{equation}
where $\mu_H$ is the harmonic mean given by $\mu_H\{a_1, \ldots, a_N  \}\, = \, \big\{\frac{1}{N} \sum_{i=1}^N a_i^{-1}   \big\}^{-1}$ for positive real numbers  $a_1, \ldots, a_N$, and where $H_m$ is the $m^{th}$ harmonic number as before.

\item[ (b)]  For $\alpha \in (0,1)$, we have the limiting expression
\begin{equation}
\label{expectationlimiting}
\lim_{n \to \infty} \mathbb{E}\big(p_2|Y=\,\lfloor \alpha n \rfloor \big) \, = \, \frac{1-\alpha}{- \ln (\alpha)}\,,
\end{equation}
where $\lfloor a \rfloor$ is the integer part of $a>0$.
\end{enumerate}
{\bf Proof.}
\begin{enumerate}
\item[ (a)]   Expression (\ref{expectationharmonic1}) follows directly from (\ref{equationbetamixture}) as  $\mathbb{E}(p_2|y) \, = \,  \frac{\sum_{k=0}^{n-y}  \frac{1}{k+y+1} \, \frac{k+y+1}{n+2}}{H_{n+1}-H_y}$, while (\ref{expectationharmonic2}) follows from part (\ref{expectationharmonic1}) since $\mu_H\{y+1, y+2, \ldots, n+1  \}\,=\,  \frac{n-y+1}{H_{n+1}-H_y}$.

\item[ (b)]  Since $\lim_{m \to \infty} \big\{H_m - \ln (m)\big\} = \gamma$, we obtain from (\ref{expectationharmonic1})
\begin{align*}
\lim_{n \to \infty} \mathbb{E}(p_2|Y\, = \lfloor \alpha n\rfloor  \,) \, & = \, \frac{\lim_{n \to \infty} (n-\lfloor \alpha n\rfloor +1)/(n+2)}{\lim_{n \to \infty} \big(H_{n+1} - H_{\lfloor \alpha n\rfloor }\big)} \, \\
\,& =  \lim_{n \to \infty} \big\{\frac{1-\alpha}{\ln(n+1) - \ln(\lfloor \alpha n\rfloor) }\big\} =\frac{1-\alpha}{- \ln (\alpha)}\,.  \qed 
\end{align*}
\end{enumerate}
\end{corollary}

\noindent  Part (b) of Corollary \ref{corollaryposteriorexpectation} does not cover the cases $y=0$ or $y=n$.  However, these are directly available from (\ref{expectationharmonic1}) as 
\begin{equation}
\nonumber
\label{posteriorexpectationy0y=n}
\mathbb{E}(p_2|0)\,=\,  \frac{n+1}{n+2} \, \frac{1}{H_{n+1}}  \, \hbox{ and } \mathbb{E}(p_2|n)\,=\,\frac{n+1}{n+2}\,.
\end{equation}
The latter of these was already available from the posterior distribution result $p_2|n \sim \hbox{Beta}(n+1,1)$ appearing in Remark \ref{remarkY=n}, while the former can be well approximated or bounded with properties of $H_{n+1}$.  For instance by making use of (\ref{harmonicapprox}), we obtain immediately the accurate approximation at least for moderate or large $n$, 
\begin{equation}
\label{posteriorexpbothy=0y=n}
\mathbb{E}(p_2|0)\, \approx \, \frac{n+1}{n+2} \, \frac{1}{\ln (n+1) + \gamma + \frac{1}{2(n+1)}} \,.
%  \hbox{ and } 
% \frac{n+1}{n+2} \, \frac{1}{\ln (n+1) + \gamma + \frac{1}{2(n+1)}} \leq \mathbb{E}(p_2|0) 
%\leq \frac{n+1}{n+2} \, \frac{1}{\ln (n+1) + \gamma + \frac{1}{2 n}}\,.
 \end{equation}
 
\noindent   Expression (\ref{expectationlimiting}) can be interpreted as a convenient large sample approximation for the posterior expectation $\mathbb{E}\big(p_2|Y=\,\lfloor \alpha n \rfloor \big)$, but (\ref{expectationharmonic1}) remains exact and leads to more accurate approximations or sharp bounds in conjunction with (\ref{harmonicapprox}) and (\ref{harmonicbounds}), respectively. For instance, for $y$ neither too small or too large, we can plug-in either 
\begin{equation}
\nonumber
H_{n+1}-H_y \simeq \ln \left(\frac{n+1}{y}\right)-\frac{n-y+1}{2y(n+1)}\,, \hbox{ or } H_{n+1}-H_y \simeq \ln \left(\frac{n+1}{y}\right)-\frac{n-y+1}{2y(n+1)}+\frac{(n+1)^2-y^2}{12(n+1)^2y^2}\,,
\end{equation}
into (\ref{expectationharmonic1}) to obtain accurate approximations for $\mathbb{E}(p_2|y)$. Similarly, the bounds
\begin{equation}
\label{bounds}
\frac{n-y+1}{n+2} \, \frac{1}{\ln (\frac{n+1}{y} ) + \frac{1}{2 (n+1)} - \frac{1}{2(y+1)}} \leq \mathbb{E}(p_2|y) \leq \frac{n-y+1}{n+2} \, \frac{1}{\ln (\frac{n+1}{y} ) + \frac{1}{2(n+2)} - \frac{1}{2y}}\,
\end{equation}
follow for $0 < y < n$.   Finally, adapting the above in an asymptotic context, we obtain the following.

\begin{corollary}
\label{corollaryposteriorexpapprox}
Consider model (\ref{model}), fixed $\alpha \in (0,1)$, and suppose that $Y=\lfloor \alpha n \rfloor$ has been observed.   Then, we have
\begin{equation}
\label{expectationlimiting2}
\mathbb{E}\big(p_2|Y=\,\lfloor \alpha n \rfloor \big) \, \approx \, \frac{(n(1-\alpha)+1)/(n+2)}{\ln(1+1/n)-\ln(\alpha) - \big(n(1-\alpha)+1\big)/\big(2n(n+1)\, \alpha  \big)} \,=\, 
\hat{E}(\alpha,n) \, \hbox{ (say) }\,,
\end{equation}
where $\lim_{n \to \infty} \frac{\mathbb{E}\big(p_2\big|Y=\,\lfloor \alpha n \rfloor \big)}{\hat{E}(\alpha,n)} \,=\, 1 $.
\end{corollary}
\noindent {\bf Proof.}
\noindent The result follows by approximating the denominator in (\ref{expectationharmonic1})
with (\ref{harmonicapprox}).  \qed

We conclude this section with some numerical evaluations.    Large sample posterior expectation approximations can be derived from (\ref{expectationlimiting}) or  (\ref{expectationlimiting2}) in a straightforward manner.  For example, with $n= 100 000$ and $y= 50 000$ as evoked by Surjanovic et al. (2025), or more generally with $Y= \frac{n}{2}$, we have $\mathbb{E}(p_2|Y=\frac{n}{2}) \approx  \frac{1}{2 \ln 2} \,=\,0.721348...$ for large $n$.   The bounds in (\ref{bounds}) are quite tight and yield: $ 0.7213397155 \leq \mathbb{E}(p_2|Y=50000) \leq 0.7213709364$ for $n=100000$.

\noindent For smaller, or much smaller, $n$, the above approximation may be rather crude and (\ref{expectationlimiting2}) is more precise.  As an illustration, Figure 1 compares 
the approximations $\frac{(1-\alpha)}{-\ln(\alpha)}$ and $\hat{E}(\alpha,n)$.  These are extended here for illustrative purpose to cases where $\alpha n$ is not an integer.  In (a),  $\alpha$ is set as $1/2$ and the approximations are contrasted as a function of $n$, while part (b) fixes $n=100$ and the approximations are comparable as a function of $\alpha$.  Clearly, the additional correction represented by $\hat{E}(\alpha,n)$ has a significant impact, at least for such small or not so large values of $n$.    Other choices of $\alpha$ in (a) and $n$ in (b), which we experimented with,  paint the same picture.

\begin{figure}[htbp]
\centering
\begin{subfigure}{0.48\textwidth}
  \centering
  \includegraphics[width=\linewidth]{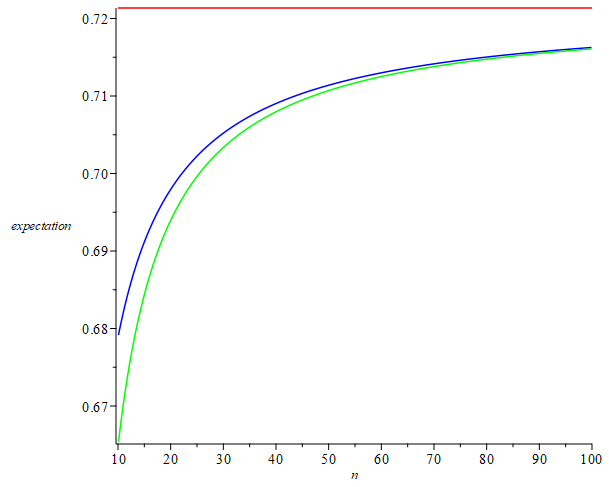}
  \caption{Expectation $\mathbb{E}(p_2|y)$ for $y=n/2$ as a function of $n$ (blue) in comparison to approximations (\ref{expectationlimiting}) (red) and (\ref{expectationlimiting2}) (green). }
  \label{fig:first}
\end{subfigure}
\hfill
\begin{subfigure}{0.48\textwidth}
  \centering
  \includegraphics[width=\linewidth]{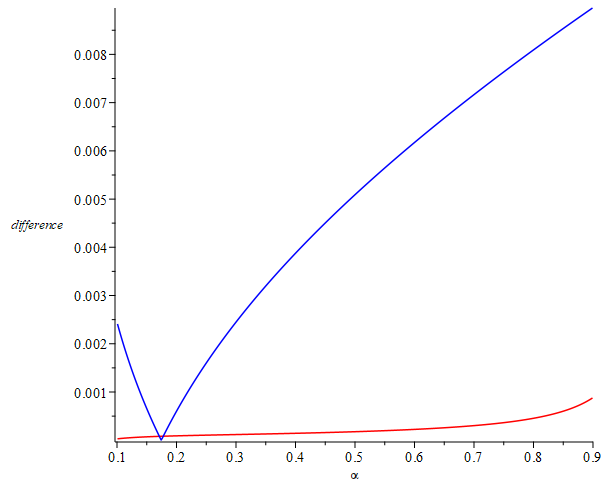}
  \caption{Differences $|\mathbb{E}(p_2|y)- \hat{E}(\alpha,100)|$ (red) and $|\mathbb{E}(p_2|y)- \frac{1-\alpha}{-\ln(\alpha)}|$ (blue), for $y=100 \alpha$ as functions of $\alpha$.}
  \label{fig:second}
\end{subfigure}
\caption{}
\label{fig:both}
\end{figure}

\section{Alternative parametrization and marginal posterior distributions}
\label{sectionalternative}

An alternative and straightforward approach is to reparametrize with $\theta_1=p_2$ and $\theta_2=p_1 p_2$. In this case, the likelihood depends on $\theta_2$ only, while the prior density impacts the posterior distributions of $\theta_2|y$ and $\theta_1|\theta_2,y$.   From (\ref{jointposteriordensity}), one obtains the decomposition $\pi(\theta_1, \theta_2|y) \, = \,  \pi(\theta_1|\theta_2, y) \, \pi(\theta_2|y) $ with 
 \begin{equation}
\label{posteriorconditional}  \pi(\theta_1|\theta_2, y) = - \, \frac{1}{\theta_1 \ln (\theta_2) }\, \mathbb{I}_{(\theta_2,1)}(\theta_1),
\end{equation}
independently of $y$, and 
\begin{equation}
\label{posteriormarginal}
\pi(\theta_2|y) \,= \, \frac{\binom{n}{y} \,\theta_2^y \, (1-\theta_2)^{n-y} \,\, (-\ln \theta_2)}{p(y)} \, \, \mathbb{I}_{(0,1)}(\theta_2)\,.
\end{equation}
In terms of the posterior expectation, observe that $\mathbb{E}(\theta_1|\theta_2, y) \, = \, \frac{\theta_2 -1}{\ln \theta_2}$ which leads to
\begin{eqnarray}
%\label{equationsecondexpression}
\nonumber \mathbb{E}(p_2|y) \, = \, \mathbb{E}(\theta_1|y) \, & = &\, \mathbb{E} \big(\frac{\theta_2 -1}{\ln \theta_2}  \, \big| y \big) \\
\nonumber & = &\, \binom{n}{y} \, \,\frac{ \int_0^1   \theta_2^y \, (1-\theta_2)^{n-y+1}  \, d\theta_2}{p(y)} \\
\label{expectationpassagetopy} & = &\, \frac{n-y+1}{(n+1) \, (n+2) \, p(y)}\,,
\end{eqnarray}
which is equivalent to (\ref{expectationharmonic1}). 

\begin{remark}
\label{remarkSibuya}
In Sibuya (1979), the posterior density $\pi(\theta_2|y)$ appears as a ``(zero)-end accented'' Beta distribution, and the normalization factor $p(y)$ is expressed as the difference of two digamma functions which is equivalent to expression (\ref{marginalpmf}).
\end{remark}

We now provide another convenient expression for the posterior density of $p_2$ (equivalently of $p_1$ or $\theta_1$), which differs from the mixture representations given in parts (a)
 and (b) of Theorem \ref{theorem}.  We furthermore establish that the posterior variance, in the asymptotic setting with $\frac{Y}{n} \to \alpha \in (0,1)$ as $n \to \infty$, converges to a positive constant, indicative of the absence of posterior contraction.  The latter is itself  
a consequence of the property that the joint posterior distribution of $(p_1,p_2)$ concentrates, under the above asymptotic setting, around the curve  $\{(p_1,p_2) \in (0,1)^2\,:\,   p_1\,p_2=\alpha \}$ rather than a single point.

\begin{theorem}
\label{theoremmarginalposterior+variance}
Consider model (\ref{model}) and suppose that $Y=y$ has been observed.  Then,
\begin{enumerate}
\item[ (a)]   The posterior density of $p_2$ (or $\theta_1$) is given by
\begin{equation}
\label{posteriordensitysection5}   \pi(p_2|y) \, = \, \frac{F_{y+1,n-y+1}(p_2)/p_2}{\, H_{n+1} - H_y} \,
,
\end{equation}
where $F_{a,b}$ is the cdf of a $\mathbb{B}(a,b,1)$ random variable;
\item[ (b)]  For $\alpha \in (0,1)$, we have the limiting expression
\begin{equation}
\label{variancelimiting}
\lim_{n \to \infty} \mathbb{V}\big(p_2|Y=\,\lfloor \alpha n \rfloor \big) \, = \, \frac{(1-\alpha)}{- \ln (\alpha)}\, \Big(\frac{1+\alpha}{2} \, - \, \frac{1-\alpha}{- \ln (\alpha)}     \Big) \,>0\,.
\end{equation}
\end{enumerate}
\end{theorem}
\noindent {\bf Proof.}  
(a)  From (\ref{posteriorconditional}), (\ref{posteriormarginal}), and (\ref{marginalpmf}), we have
\begin{align*}
\pi(\theta_1, \theta_2) \, & = \, \frac{1}{p(y)} \, \frac{1}{\theta_1} \, \binom{n}{y} \, \theta_2^y (1-\theta_2)^{n-y} \, \mathbb{I}_{(0,\theta_1)}(\theta_2) \, \mathbb{I}_{(0,1)}(\theta_1) \\
&  =  \frac{1}{(n+1) \, p(y) \, \theta_1} \,  f_{y+1,n-y+1,1}(\theta_2) \\
& = \frac{1}{(H_{n+1}\,-\,H_y) \, \theta_1} \,  f_{y+1,n-y+1,1}(\theta_2)\,.
\end{align*}
The result then follows by integrating out $\theta_2$. \\
\noindent (b)  Direct calculations with the density in part (a) and the second moment of a Beta random variable yield
\begin{align*}
\mathbb{E}(p_2^2|y) \, = \, \mathbb{E}(\theta_1^2|y) \, &  = \, \frac{1}{(H_{n+1}\,-\,H_y)} \, \int_{0}^{1} \int_{\theta_2}^{1} \theta_1   f_{y+1,n-y+1,1}(\theta_2) \, d\theta_1 \, d\theta_2 \\
 \, &  = \, \frac{1}{(H_{n+1}\,-\,H_y)} \, \int_{0}^{1} \, \frac{1-\theta_2^2}{2} \, f_{y+1,n-y+1,1}(\theta_2) \, d\theta_2 \\
 \, &  = \, \frac{1}{2(H_{n+1}\,-\,H_y)} \, \Big(1 \, - \, \frac{(y+1)(y+2)}{(n+2)(n+3)} \Big)\,.
\end{align*}
With (\ref{harmoniclimits}), we obtain $\lim_{n \to \infty} \mathbb{E}\big(p_2^2|Y=\,\lfloor \alpha n \rfloor \big) \, = \, \frac{1-\alpha^2}{-2 \ln(\alpha)}$, and the result then follows with part (b) of Corollary \ref{corollaryposteriorexpectation} which tells us that $
\mathbb{E}^2\big(p_2|Y=\,\lfloor \alpha n \rfloor \big) \, = \, \frac{(1-\alpha)^2}{ \,\ln^2(\alpha)}$. \qed

\begin{remark}
\label{remarky=0,y=n}   The cases $\alpha=0$ and $\alpha=n$ (i.e., $Y=0$ and $Y=n$) are not covered by (\ref{variancelimiting}), but they are handled directly, for instance via (\ref{posteriordensitysection5}). In both cases, the marginal posterior distribution of $p_1$ concentrates asymptotically towards a point, respectively $0$ and $1$ depending on whether $Y=0$ or $Y=n$.   Indeed, as already seen in Remark \ref{remarkY=n}, $p_2|n \sim \mathbb{B}(n+1,1)$ with $\mathbb{E}(p_2|n)\,=\, \frac{n+1}{n+2} \to 1 $ and $\mathbb{V}(p_2|n)\,=\, \frac{n+1}{(n+2)^2 (n+3)} \to 0$ as $n \to \infty$.   For $Y=0$, we have from above (\ref{posteriorexpbothy=0y=n}):  $0 \leq \mathbb{V}(p_2|0) \leq \mathbb{E}(p_2^2|0) \leq 
\mathbb{E}(p_2|0) = \frac{n+1}{n+2} \, \frac{1}{H_{n+1}} \to 0$ as $n \to \infty$.
\end{remark}

\begin{remark}
\label{approximationforp(y)}
With the normalization constant $p(y)$ and posterior expectation $\mathbb{E}(p_2|y)$ linked by
equation (\ref{expectationpassagetopy}), so are various approximations and bounds.  For instance, part (b) of Corollary \ref{corollaryposteriorexpectation} implies that
\begin{equation}
p(y) \approx \frac{n(1-\alpha)+1}{(n+1) \, (n+2) \, (1-\alpha)/-\ln(\alpha)} \approx  \frac{-\ln(\alpha)}{n}
\end{equation} 
for $y=\alpha n$ and large $n$, while a better approximation is attainable with Corollary \ref{corollaryposteriorexpapprox} as:
\begin{equation}
\nonumber   p(y) \approx \frac{n(1-\alpha)+1}{(n+1) \, (n+2) \, \hat{E}(\alpha,n)}\,. 
\end{equation}
Revisiting a numerical simulation provided by Surjanovic et al. (2025) where 
$y=50 000$ and $n=100 000$, the first of the above equations yields   $\ln p(50 000) \, \approx \, \ln\big\{(\ln 2) \times 10^{-5} \big\} \, = \, -11.879438..$ which corroborates Surjanovic et al. (2025)'s given value $\ln p(50 000) \approx -11.8794...$.
\end{remark}

We conclude by discussing large sample approximations for the posterior density $\pi(p_2|y)$.  With the cdf of a Beta distribution approximated by a normal distribution cdf with matching mean and standard deviation, we have for $y=\alpha n, \alpha \in (0,1),$
\begin{equation}
\nonumber
F_{y+1,n-y+1}(p_1) \, \approx \, \Phi_{\alpha,\sqrt{\alpha(1-\alpha)/n}}\,(p_1)
\end{equation}
where $\Phi_{\mu,\sigma}(\cdot)$ is the $N(\mu,\sigma^2)$ cdf.   From the above, part (a) of Theorem \ref{theoremmarginalposterior+variance}, and the first approximation displayed in (\ref{harmonicapprox}), we obtain
\begin{equation}
\label{posteriorapproximate2}
\pi(p_2|y) \approx \frac{1}{C p_2}\Phi_{\alpha,\sqrt{\alpha(1-\alpha)/n}}\,(p_2)\,,
\end{equation}
with 
$$C \, \approx \, H_{n+1}-H_y \, \approx \, \ln(1+1/n)-\ln(\alpha) - \big(n(1-\alpha)+1\big)/\big(2n(n+1)\, \alpha  \big)\, = \, \hat{C}\, (\hbox{ say })\,.$$
Here $C$ is a normalization constant.  Figure 2 compares for $(n,y)=(20,10)$ the exact posterior density with the approximation $\hat{\pi}(p_2|y)\,=\, \frac{1}{\hat{C}}
  \frac{1}{p_2} \, \Phi_{\alpha,\sqrt{\alpha(1-\alpha)/n}}\,(p_2)$  obtained by substituting $\hat{C}$ for $C$ in 
(\ref{posteriorapproximate2}).  The latter is not quite a density, but very close here as supported by a numerical evaluation which yields an area $\in (0.9999, 1.0000)$.  Given that $n=20$ only, the proximity of the two curves is striking.   The integrated $L_1$ distance $\int_0^1 |\pi(p_2|y) \, - \, \hat{\pi}(p_2|y)| \, dp_2$ is about equal to $0.0219$.  To conclude, we point out  that the density in (\ref{posteriorapproximate2}) converges as $n \to \infty$ to the density
\begin{equation}
\nonumber  \pi^*(p_2|y) \, =\, \frac{1}{- p_2 \, \ln(\alpha)} \, \mathbb{I}_{(\alpha,1)}(p_2)\,.
\end{equation}
From this, it is interesting to observe that: (i) the matching of the expectation and variance associated with $\pi^*$ with the asymptotic expressions
(\ref{expectationlimiting}) and (\ref{variancelimiting}), (ii) the presence of the limiting support $(\alpha,1)$, and (iii) (again) the absence of posterior contraction.  
\begin{figure}[htbp]
\centering
%\begin{subfigure}{0.48\textwidth}
  \centering
  \includegraphics[width=0.66\textwidth]{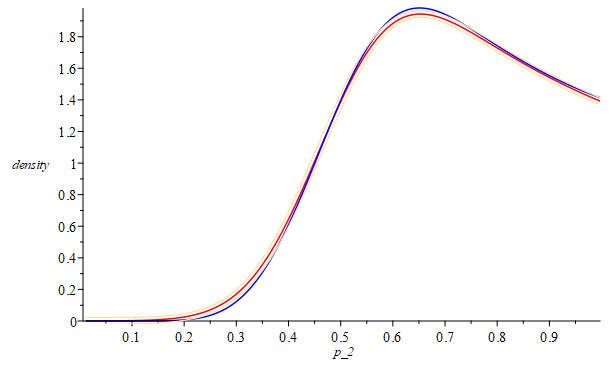}
  \caption{Posterior density $\pi(p_2|y)$ (blue) and the approximation $\frac{1}{\hat{C}}
  \frac{1}{p_2} \, \Phi_{\alpha,\sqrt{\alpha(1-\alpha)/n}}(p_2)$ (red) for $n=20, y=10$.    }
  \label{fig:third}
\end{figure}
\section*{Acknowledgements}
The authors wish to thank two reviewers for precise and pertinent comments that led to an improved presentation of the findings.   The authors are grateful to Alexandre Bouchard-C\^ot\'e and Jim Berger for useful exchanges. Author Marchand acknowledges and is thankful for research support from NSERC of Canada.

\section*{Appendix}

Here is an extension of Lemma \ref{lemmamarginalYchris} which involves a conditional binomial model with $l$ uniformly distributed probabilities $p_1, \dots, p_l$.  The result applies for all $l \geq 2$ and yields an alternate expression for $l=2$.
\begin{lemma}
\label{lemmamarginalYgeq2}
  Consider the following extension of model (\ref{model}):
\begin{equation}
\label{modelextension}
Y|p_1,\dots, p_m \sim \hbox{Binomial}(n, \prod_{i=1}^m p_i) \;, \hbox{with } \;  p_1, \dots, p_l \sim^{\!\!\!\!\!\!\!\! i.i.d.} \mathcal{U}(0,1)\,, 
\end{equation}  
\end{lemma}  
  The marginal probability mass function of $Y$ is given by
\begin{equation}
\label{marginalpmfseveral}
p(y) \, = \, \frac{n!}{y!} \, \sum_{k=0}^{n-y}  \frac{(-1)^{k}}{k! (n-y-k)!} \, \frac{1}{(y+1+k)^l}\, \mathbb{I}_{\{0,1, \ldots, n\}}(y)\,.
\end{equation}  
{\bf Proof.} With $W=\prod_{i=1}^l p_i$ having density $\frac{(-\ln(w))^{l-1}}{\Gamma(l)} \, \mathbb{I}_{(0,1)}(w)$ and $Y|W \sim \hbox{Binomial}(n, W)$, we have 
\begin{eqnarray*}
p(y) \, & = & \, \frac{1}{\Gamma(l)} \, \int_0^1  \binom{n}{y} \, w^y (1-w)^{n-y} \, (-\ln(w))^{l-1} \, dw \\
\, & = & \, \frac{1}{\Gamma(l)} \, \, \binom{n}{y} \, \int_0^{\infty} z^{l-1} \, e^{-z(y+1)} \, (1-e^{-z})^{n-y} \, dz\\
\, & = & \, \frac{1}{\Gamma(l)} \, \binom{n}{y} \, \int_0^{\infty} z^{l-1} \, e^{-z(y+1)} \, \sum_{k=0}^{n-y} \, \binom{n-y}{k} \, (-1)^k \, e^{-kz}  \, dz\\
\, & = &  \frac{1}{\Gamma(l)} \, \frac{n!}{y!} \,\, \sum_{k=0}^{n-y} \, \frac{(-1)^{k} }{k! (n-y-k)!} \, \int_0^{\infty}  z^{l-1} \, e^{-z(y+1+k)}  \, dz \,,
\end{eqnarray*}
which yields the result. \qed

\noindent  Observe that $p(n)\,=\, \frac{1}{(n+1)^l}$ extends nicely the earlier result for $l=2$.    The above provides as well an alternate expression for the marginal distribution of $Y$ when $l=2$, and a direct passage from (\ref{marginalpmfseveral}) to (\ref{marginalpmf}) is available through a partial fraction decomposition.
\end{document}